\documentclass{amsart}
\usepackage{amsfonts,amscd}
 
\newtheorem{claim}{}[section]
\newtheorem{theorem}[claim]{Theorem}
\newtheorem{lemma}[claim]{Lemma}
\newtheorem{proposition}[claim]{Proposition}
\newtheorem{corollary}[claim]{Corollary}

\newtheorem{definition}[claim]{Definition}

\def\proclaim #1. #2\par{\medbreak
\noindent{\bf#1.\enspace}{\sl#2}\par\medbreak}
\makeatother

\DeclareMathOperator{\Bdb}{{\mathbb B}}
\DeclareMathOperator{\Kdb}{{\mathbb K}}
\DeclareMathOperator{\Cdb}{{\mathbb C}}
\DeclareMathOperator{\Rdb}{{\mathbb R}}

\DeclareMathOperator{\Ndb}{{\mathbb N}}

\DeclareMathOperator{\A}{{\mathcal A}}
\DeclareMathOperator{\M}{{\mathcal M}}
\DeclareMathOperator{\Sy}{{\mathcal S}}

\begin{document}
 
\title[Duality and operator algebras]{Duality and operator
 algebras}
 
\date{July 16, 2004}
 
\author{David P. Blecher}
\address{Department of Mathematics, University of Houston, Houston,
TX
77204-3008}
\email[David P. Blecher]{dblecher@math.uh.edu}
 \author{Bojan Magajna} 
\address{Department of Mathematics,
University of Ljubljana,
Jadranska 19, Ljubljana 1000, Slovenia}
\email[Bojan Magajna]{Bojan.Magajna@fmf.uni-lj.si}
\thanks{*Blecher is partially supported by a grant from
the National Science Foundation.
Magajna is  partially supported by the Ministry of Science and
Education of Slovenia. } 

\begin{abstract}
We investigate some subtle and interesting
phenomena in the duality theory
of operator spaces and operator algebras.  
In particular,
we give several applications of operator space
theory, based on the surprising fact that 
certain maps are always weak$^*$-continuous on 
dual operator spaces.    For example, 
 if $X$ is a subspace of a $C^*$-algebra $A$, and if
$a \in A$ satisfies $a X \subset X$ and $a^* X \subset X$,
and if $X$ is isometric to a dual 
Banach space, then we show that the function $x \mapsto a x$ on $X$
is weak$^*$ continuous.  
Applications include a new
characterization of the $\sigma$-weakly closed (possibly nonunital and
nonselfadjoint) operator algebras,
and it makes possible a generalization of the theory of $W^*$-modules  
to the framework of modules over such algebras.
We also give a Banach module characterization of
$\sigma$-weakly closed spaces of operators which are
invariant under the action of a von Neumann algebra.
\end{abstract}
 
\maketitle

\let\text=\mbox
 
\section{Introduction}

Functional analytic questions about spaces of operators often 
boil down to considerations involving dual, or weak$^*$, topologies.
In many such calculations, the key point is to prove that 
certain linear functions are weak$^*$ continuous.   In 
the present paper we offer a couple of results which ensure that  
a linear map be automatically continuous with respect to
such topologies.  In the right situation, these results can   
be extremely useful.
  
A {\em multiplier} of a
Banach space $E$ is a linear map $T : E \to E$ such that 
there exists an isometric embedding $\sigma :
E \rightarrow C(\Omega)$
for a compact space $\Omega$, and a function $a \in C(\Omega)$,
such that $\sigma(T x) = a \sigma(x)$ for all $x \in E$ 
(see \cite[Section I.3]{MIBS} and \cite[Theorem 3.7.2]{BLM}).  If, 
further, $E$ is a dual Banach space, we are able to show the surprising
fact that $T$ is automatically $w^*$-continuous.  
It is clear how to generalize the notion
of multipliers, to maps on an operator space $X$.
By definition, an operator  space is a subspace of a $C^*$-algebra.
Thus, a {\em left multiplier} of an operator space $X$ is
a linear map $T : X \to X$ such that
there exists a {\em completely isometric} (this term is
defined below) embedding $\sigma : X \rightarrow A$,
for a $C^*$-algebra $A$, and an element $a \in A$,
such that $\sigma(T x) = a \sigma(x)$ for all $x \in X$.
These operator space multipliers have been useful
in various ways in the last several years (see
e.g.\ \cite{MCAA,TCOO} or \cite[Chapter 4]{BLM} and
references therein).
 We shall see here that left multipliers 
of dual operator spaces
are automatically $w^*$-continuous too, and we shall 
give several remarkable applications of this fact,
mostly to {\em operator algebras} (that is, 
possibly nonselfadjoint subalgebras of $C^*$-algebras). 
We are also able to relax the restriction above that 
$\sigma$ be a complete isometry, 
and allow $X$ to be a
dual Banach space and $\sigma$ an isometry,
 provided that $a^* \sigma(X) \subset \sigma(X)$ too.
This    has some interesting consequences.
   
Generally, our paper is concerned with some subtle 
and interesting phenomena in the duality theory of
operator spaces and operator algebras.  For example,
in Section 2 we give perhaps the simplest example
of an operator space which is a dual Banach space but
not a dual operator space.  
In Section 3, we establish a property of projection maps
 on an operator module.  This we use as a
key technical tool, although it is of interest in its own right. 
We also prove our first `automatic continuity' result.
 In Section 4, we prove our main result, that left multipliers
of dual operator spaces are automatically $w^*$-continuous.
This is surprising in the light of any of the known 
alternative definitions of left multipliers (see \cite[Theorem 4.5.2]{BLM}
or \cite{MCAA}).
 We also give numerous corollaries and complementary results.   
For example, we give a new characterization of 
$\sigma$-weakly closed algebras of operators,
and also a Banach module characterization of
$\sigma$-weakly closed spaces of operators which are
invariant under the action of a von Neumann algebra.
Finally, in Section 5, we illustrate again the power
of our main result, by using it to generalize key aspects
of the theory of the `selfdual
modules' of Paschke \cite{IPMO} and Rieffel
\cite{MEFC}, also known as
$W^*$-modules, to nonselfadjoint operator algebras.       
 
We now turn to the definitions, and background facts.
The basic source text used for background information
is \cite{BLM}, which should be available soon (it has a
2004 publication date).
Some of this information may be found in \cite{CBMAO} too.
Because \cite{BLM} is temporarily unavailable to all, 
in this preprint version of the present paper
we will also reference background results using 
\cite{CBMAO} too, where possible.
In any case, see e.g.\ these two
texts for more explanation of the notation and
facts below, if needed.

Throughout $H$ and $K$ are Hilbert spaces.
Operator spaces, defined above, may also be 
thought of as the closed linear subspaces of $B(K,H)$.
Equivalently, by  a theorem of Ruan, an operator space
is a vector space 
$X$ with a norm defined on each of
the spaces $M_n(X)$ of $n \times n$ matrices over 
$X$ satisfying two compatibility conditions
which we shall not spell out here.  
A linear map $T : X \to Y$ between operator spaces clearly 
induces a map $T_n : M_n(X) \to M_n(Y)$.  We say that 
$T$ is {\em completely 
isometric} (resp.\ {\em completely contractive},
{\em a 
complete quotient map}) if $T_n$ is
isometric (resp.\ contractive, takes the open ball of $M_n(X)$ onto the open
ball of $M_n(Y)$), for all $n \in \Ndb$.
We say that
$T$ is {\em completely  bounded} if 
$$
\Vert T \Vert_{\rm cb} \; \overset{def}{=} \;  \sup_{n \in \Ndb}  \Vert T_n \Vert
 \; < \; \infty.
$$
We write $CB(X,Y)$ for the space of
completely  bounded maps, with this norm,
and $CB(X) = CB(X,X)$.  If $X$ and $Y$ are left (resp.\ right)
$A$-modules
then $_ACB(X,Y)$ (resp.\ $CB_A(X,Y)$)
denotes the subspace of $CB(X,Y)$ consisting
of module maps.  The reader can guess the meaning
of $B_A(X,Y)$, $_ACB(X)$, etc.
   We say that an operator space $X$ is {\em unital} if it has a
distinguished element $1$,
such that there exists a complete isometry $T \colon X
\rightarrow A$ into a unital $C^*$-algebra with $T(1) = 1_A$.
Examples of these 
include {\em operator systems}, namely linear selfadjoint subspaces
of a $C^*$-algebra $A$ with $1_A \in X$.
    
Operator algebras (defined in the second paragraph of our 
paper) may also be 
defined purely abstractly in terms of matrix norms
(e.g.\ see \cite[Theorem 2.3.2]{BLM} or 
\cite[p.\ 252]{CBMAO}).  We say that an 
operator algebra is {\em approximately unital}
if it has a contractive approximate identity.
  If  $X$ is an operator space, then the set $\M_l(X)$
of left multipliers (also defined in the second paragraph) of $X$,
turns out to be such an abstract operator algebra.
There are several equivalent
definitions of  $\M_l(X)$ in the literature (e.g.\ see \cite{MCAA}
or \cite[Sections 4.5 and 8.4]{BLM}
or \cite[Chapter 16]{CBMAO}).  
For example, for a unital operator space $X$, one may define 
$\M_l(X)$ to be the image in $CB(X)$ of the subalgebra
$\{ a \in D : a X \subset X \}$, via the canonical map
taking such $a$ to the map $x \mapsto ax$ on $X$.  Here 
$D$ is a certain `extremal' $C^*$-algebra containing $X$, with
$1 = 1_D$.  In fact, $D$ may be taken to be either 
Hamana's {\em injective envelope}, or {\em $C^*$-envelope} 
(also known as Arveson's {\em noncommutative Shilov boundary}),
of $X$.  See e.g.\ \cite{SOC,IEOO},  
\cite[Sections 4.2--4.4, and 8.3]{BLM}
or \cite[Chapter 15]{CBMAO}
 for a thorough discussion of of the latter objects.

A dual Banach space is a Banach space linearly
isometric to the dual of another Banach space (the latter is
called a predual).   We abbreviate the word
`weak*' to `$w^*$'.  The $w^*$-topology on $B(H)$ is
often called the $\sigma$-weak topology.   The 
product of $B(H)$ (and hence of any
$w^*$-closed subalgebra of $B(H)$)
is a separately $w^*$-continuous bilinear map.     
A $W^*$-algebra is a $C^*$-algebra which is
a dual Banach space;  by a 
well known theorem of Sakai, these are `exactly'
the von Neumann algebras.
Indeed, the methods of our paper owe enormously to Tomiyama's quick proof
of Sakai's theorem, and adaptions of this
method by others, e.g.\ 
\cite[Theorem 9.1]{Strat}, \cite{OIAN}.
The second dual of a $C^*$-algebra is a $W^*$-algebra
(a fact for which there exist simple proofs in the literature).
A consequence of the well-known Krein-Smulian theorem, is that 
a linear bounded map $u\colon E \to F$ between
dual Banach spaces is $w^*$-continuous if and only if whenever
$x_t \rightarrow x$ is a bounded net converging in the
$w^*$-topology in the domain space, then $u(x_t) \rightarrow u(x)$
in the $w^*$-topology.  If this holds, 
and if moreover $u$ is a $w^*$-continuous
isometry, then $u$ has $w^*$-closed range, and $u$ is a
$w^*$-$w^*$-homeomorphism onto {\rm Ran}$(u)$.   See
\cite{ROOB} or  \cite{MADO}
for proofs.  These facts will be used silently very often in our paper.

If $A$ is a $C^*$-algebra, we write $M(A)$ for its
multiplier algebra.
If $X$ and $Y$ are sets (in a $C^*$-algebra say)
then we write 
$X Y$ for the {\em norm closure} of the span of terms 
of the form $x y$, for $x \in X, y \in Y$.   Similar conventions
hold for products of three subsets.

\section{Dual operator spaces}

The Banach space dual $Y^*$ of an operator space $Y$ is again an 
operator space, in a canonical way.  Namely, for $n \geq 2$
we assign $M_n(X^*)$ the norm pulled back via the canonical
algebraic isomorphism $M_n(X^*) \cong CB(X,M_n)$  (e.g.\ see 
\cite[Section 1.4]{BLM} or \cite{MADO}).   
We recall that $X$ is a {\em dual operator space} if $X$ is
completely isometric to $Y^*$, for an operator space $Y$.
Le Merdy gave a beautiful characterization of 
dual operator spaces (see \cite{OTDO} and 
1.6.4 in \cite{BLM});
and he also showed that an operator space   
which is a dual Banach space need not be a dual operator space.
Simpler examples were later found by Peters-Wittstock, Effros-Ozawa-Ruan
 \cite{OIAN}, and in \cite[Remark 7.10]{DANP}.
This phenomenon will play an important role in this paper,
for example we will often ask when results valid for
dual operator spaces are also valid for an operator space
which is a dual Banach space.
Indeed, the following example, which will play a role later in
the paper, may be the simplest example of this  
phenomenon:

\begin{proposition} \label{nonp}  There is an operator space
structure on $B(\ell_2)$, for which there exists
a predual Banach space, but not a predual operator space.
\end{proposition}

\begin{proof}  Let $H = \ell_2$, let $S^\infty$ denote the 
compact operators on $\ell_2$, let $Q$ be the Calkin algebra
$B(H)/S^\infty$, and let $Q^{\rm op}$ denote its `opposite 
$C^*$-algebra'.  Let $X$ be the subspace of $B(H) \oplus^\infty
Q^{\rm op}$  consisting of pairs $(x,\dot{x})$ for $x \in B(H)$,
where $\dot{x}$ is the class of $x$ in the quotient.  Then 
$X$ is a unital operator space (in fact it is even an operator system)
which is linearly isometric 
to $B(H)$.   Thus $X$ has a predual Banach space,
the trace class $S^1$, which is even a unique predual.
However $X$ is not a  dual operator space.   
The reason for this is that the canonical embedding
$\iota : S^\infty \hookrightarrow X$ is a complete isometry.
Thus if there were an  operator space structure on 
$S^1$ such that the canonical map
$(S^1)^* \to X$ was a complete isometry, 
then the unique $w^*$-continuous contraction
$B(H) \to X$ extending $\iota$, would be a 
complete contraction (see 1.4.8 in \cite{BLM}).   This unique extension 
must be the canonical `identity' map from $B(H)$ to $X$. 
 The fact that it is
completely contractive forces the 
canonical quotient map $B(H) \to Q^{\rm op}$ to be
a complete contraction, which in turn implies that
the `identity map' $Q \to Q^{\rm op}$ is a
complete contraction.  However it is well known that
the `identity map' from a $C^*$-algebra $A$ to its 
opposite algebra $A^{\rm op}$ is a complete contraction if and 
only if the $C^*$-algebra is commutative (indeed this is
clear if one applies a `noncommutative Banach-Stone theorem'
such as \cite[Corollary 1.3.10]{BLM}
to the canonical map from $A$ to $A^{\rm op}$).  
\end{proof}

In our paper we shall be quite concerned with the multiplier
algebras of a dual operator space
$X$.  As we mentioned in Section 1 for unital operator spaces
(and a similar thing is true for general operator spaces), 
$\M_l(X)$ may be defined in terms
of either the injective envelope or the $C^*$-envelope. 
If either of the
latter two objects were a $W^*$-algebra, then many of the 
technical difficulties which we will need to overcome 
in this paper, would disappear.  Unfortunately this is not generally
the case.  To show that the methods of our paper are not 
gratuitous, it seems worthwhile to take the time to
exhibit a simple explicit example of this phenomena.

\begin{proposition}  There exists a finite dimensional
unital operator algebra $M$, such that neither its
injective envelope, nor its $C^*$-envelope, are $W^*$-algebras.
\end{proposition}

\begin{proof}  Let $X$ be the span of $\{1,x,x^2 \}$ in $C([0,1])$.
This is a unital operator space, and it generates $C([0,1])$
by the Stone-Weierstrass theorem.   It is easy to see that
$[0,1]$ is the Shilov boundary for $X$ in $[0,1]$
(because for any nontrivial closed subset $C \subset
[0,1]$, there is a function in $X$ that peaks outside of
$C$).  Thus the $C^*$-envelope of $X$ is $C([0,1])$
(by e.g.\ 4.3.4 in \cite{BLM} or 
\cite[p.\ 221]{CBMAO}).  Let $D$ be an
injective envelope of $X$ as mentioned in Section 1,
thus $D$ is a unital $C^*$-algebra with $1_D = 1$.  In fact $D$  
is also the injective envelope of $C([0,1])$ (this follows
from the basic theory of the injective envelope
from e.g.\ \cite[Section 4.2 and 4.3]{BLM} or
\cite[Chapter 15]{CBMAO}, since the
$C^*$-envelope of $X$ may be defined to be the
$C^*$-subalgebra of $D$ generated by
$X$, and any minimal $C([0,1])$-projection
on $D$ is an $X$-projection and is consequently the
identity).    Next, let $M$ be the subalgebra
of $M_2(C([0,1]))$ with $0$ in the 2-1 entry, scalars
on the main diagonal, and an element from $X$ in the 1-2 entry.
This is a five dimensional unital operator algebra.  Note that
$M + M^*$ is the Paulsen system $\Sy(X)$ (see \cite[Lemma 8.1]{CBMAO} or
\cite[p. 21]{BLM}).  By 4.2.7, 4.3.6, and 4.4.13  in
\cite{BLM},  $I(M) = I(M + M^*) = I(\Sy(X)) = M_2(D)$.
Now $M \subset M_2(C([0,1]))$, and the 
latter is a $*$-subalgebra of $M_2(D)$. Thus $C^*_e(M)$, 
the $C^*$-algebra generated by
$M$ in $M_2(D)$, is also the $C^*$-algebra generated by
$M$ in $M_2(C([0,1]))$.  Hence $C^*_e(M) = M_2(C([0,1]))$.
The second assertion of the theorem  is now clear.
The first follows too, if $D$ is not a
$W^*$-algebra.  However, Claim 1:
the injective envelope of $C([0,1])$ is the
well-known Dixmier algebra.   This is the algebra $C(Y)$
in \cite[Exercise 5.7.21]{KR}),
which is shown there to not be 
a $W^*$-algebra.   Since we are not aware
of a proof of Claim 1 in the literature, we provide one
(the paper \cite{Gons} proves the same fact, but in a different
category).
 We first note that the Dixmier algebra $C(Y)$
is injective in the category
of Banach spaces (see
\cite[Exercise 5.7.20 (viii)]{KR} and
\cite[Exercise III.1.5]{Tak}).  Hence it is
injective as an operator space (by e.g.\ 4.2.11 in \cite{BLM}).
It is easy to see that the canonical injection
$C([0,1]) \to C(Y)$ is (completely) isometric,
thus we identify
$C([0,1])$ with its image under this map.  Claim 2:
every selfadjoint element $k$ in $C(Y)$ is the least upper bound
of the functions $h \in C([0,1])$ with $h \leq k$.  If Claim 2
 holds, then it is easy to see that $C(Y)$ has the `rigidity
property', and hence is the injective envelope (see \cite[Section
4.2]{BLM} or \cite[Chapter 15]{CBMAO}).
  Indeed, suppose that
 $T : C(Y) \to C(Y)$ is a (complete)
 contraction extending the identity map on
$C([0,1])$.  Since  $T(1) = 1$, $T$ is positive.
If $k, h$ are as above, then
$h = T(h) \leq T(k)$.  By the claim, $k \leq T(k)$.
Similarly, $-k \leq T(-k)$, and so $T(k) = k$.  This proves Claim 1.
  
Claim 2 is no doubt well known, but again we are not
aware of a good reference for it.  To see it, let
$f$ be the least upper bound in
$C(Y)$ of the functions $h \in C([0,1])$ with $h \leq k$.
We may assume that $k \geq 0$ (by adding a scalar multiple
of the identity, if necessary).
 Hence $k$ is the equivalence class
of a nonnegative bounded Borel function $g$, modulo
functions which are zero except on a
meager Borel set.    By basic
measure theory we may write $g = \sum_{i=1}^\infty c_i \chi_{A_i}$,
where $c_i > 0$ and $A_i$ are Borel sets in $[0,1]$.
This sum converges pointwise to $g$.
By  \cite[Exercise 5.7.20 (iii)]{KR}, there is an open set $U_i$
in $[0,1]$, and
a meager set $N_i$ such that
$\chi_{A_i} = \chi_{U_i}$ outside of $N_i$.
Let $E = \cup_i N_i$, another meager set,
then $\sum_{i=1}^\infty c_i \chi_{U_i}$ converges
pointwise outside $E$ to $g$.
For each $i$,
let $f^{i}_n$ be an increasing sequence of continuous
nonnegative functions converging to $\chi_{U_i}$.
Let $s^N_n = \sum_{i=1}^N c_i f^{i}_n$.
Outside of $E$, $s^N_n \leq g$, so that $s^N_n \leq k$ in $C(Y)$.
Hence $s^N_n \leq f$ in $C(Y)$.
If $g'$ is a nonnegative bounded Borel function
with equivalence class $f - s^N_n$,
and if we set $r = g' + s^N_n$, then the equivalence class in
$C(Y)$ of $r$ is $f - s^N_n + s^N_n = f$.   Now $r \geq s^N_n$, so 
taking the limit over $n$ and $N$ we see that
$r \geq g$ pointwise outside of $E$.
Hence  $k \leq f$ in $C(Y)$.
  \end{proof}

{\bf Remark.}  Other such examples, at least in
the operator space case, are probably known privately to experts.
For example, one may deduce such examples from the intricate
\cite[Corollary 3.8]{REOC} (here $X$ is the span of the
 generators in the reduced $C^*$-algebra of the free 
group on $n$ generators).  

\section{Modules, and a result of Zettl}

Let $A$ be an operator algebra, and $\pi : A \to B(H)$ a
completely contractive representation.   
A concrete left operator $A$-module is a subspace
$X \subset B(K)$ such that $\pi(A) X \subset X$.
An {\em (abstract)} operator $A$-module is an 
operator space $X$ which is also an $A$-module,
such that $X$ is completely
isometrically isomorphic, via an $A$-module map,
to a concrete operator $A$-module.  Note that 
in this case it is then clear from the definitions
that the map $x \mapsto ax$ on $X$ is in $\M_l(X)$
for any $a \in A$.
 There is an elegant   
characterization of operator $A$-modules, due to
Christensen--Effros--Sinclair (c.f.\ \cite[Theorem
3.3.1]{BLM}), but we shall not need this.  Many of
the most important modules over operator algebras
are operator modules, such as Hilbert $C^*$-modules.

The main trick in the following comes from 
a well known proof of Tomiyama's theorem of conditional
expectations \cite[Theorem 9.1]{Strat}.  This trick is also used 
in  \cite[Theorem 2.5]{OIAN}:
 
 \begin{lemma}   \label{newl}  Suppose that $X$ is a 
left Banach module over a $C^*$-algebra $A$, which is isometrically $A$-isomorphic 
to a nondegenerate left
 operator $A$-module.  Suppose, further,  that $Y$
is a submodule of $X$ for which there exists a contractive
linear  projection $\Phi$ from $X$ onto $Y$.    Then $\Phi$ is
an $A$-module map.
\end{lemma}  \begin{proof}   We may assume 
that $X$ is a nondegenerate left
 operator $A$-module.  We then claim that it
suffices to assume that $A$ is a $W^*$-algebra.  Indeed,
by 3.8.9 in \cite{BLM}, $X^{**}$ is a nondegenerate
left operator $A^{**}$-module, and by routine arguments $Y^{**}$
may be viewed canonically as
an $A^{**}$-submodule of $X^{**}$.  Then $\Phi^{**}$ is
a  contractive linear  projection from $X^{**}$ onto $Y^{**}$.  If
the lemma were true in the $W^*$-algebra case then
$\Phi^{**}$ is an $A^{**}$-module map, so that $\Phi$ is
an $A$-module map.
 
We next claim that it suffices to show that
\begin{equation} \label{pperp}
p^\perp \Phi(p x) \; = \; 0 , \end{equation}
for all $x \in X$, and orthogonal projections
$p \in A$.  For if (\ref{pperp}) holds then we have
$$p \Phi(x)  =  p \Phi(p x) +  p \Phi(p^\perp x) =
p \Phi(p x) = (p + p^\perp) \Phi(px) =  \Phi(px) ,$$
using (\ref{pperp}) twice (once with $p$ replaced by
$p^\perp$).  Since $A$ is densely spanned by its
projections, we conclude that $a \Phi(x) =  \Phi(ax)$
for all $a \in A$.
 
To prove (\ref{pperp}), let $y = px + t p^\perp \Phi(px)$,
where $t \in \Rdb$.  By the $C^*$-identity we have:
$$\Vert y \Vert^2 = \Vert y^* y \Vert =
\Vert x^* p x + t^2 \Phi(px)^* p^\perp \Phi(px) \Vert
\leq \Vert p x \Vert^2 + t^2 \Vert p^\perp \Phi(px) \Vert^2.$$
On the other hand, since Ran$(\Phi)$ is an $A$-submodule,
$p^\perp \Phi(px) \in {\rm Ran}(\Phi)$, so that
$$p^\perp \Phi(y) = p^\perp \Phi(px) + t p^\perp p^\perp \Phi(px)
= (1+t) p^\perp \Phi(px) .$$
Since $\Phi$ is a contraction, it follows that
$$(1+t)^2 \Vert p^\perp \Phi(px) \Vert^2 \leq
\Vert p x \Vert^2 + t^2 \Vert p^\perp \Phi(px) \Vert^2.$$
This implies that $2t \Vert p^\perp \Phi(px) \Vert^2 \leq \Vert p x \Vert^2$  
for all $t > 0$.  This is possible if and only if
$p^\perp \Phi(px) = 0$.    \end{proof}

{\bf Remark.}  The Lemma fails if $A$ is a nonselfadjoint 
operator algebra.  Indeed, if  $H$ is a Hilbert
space on which such $A$ is completely contractively represented,
and if $K$ is a closed $A$-invariant subspace of $H$,
then the orthogonal projection  onto $K$ is not necessarily an
$A$-module map.  

\medskip

\begin{theorem} \label{opwe}  Suppose that $X$ is a
left Banach module over a $C^*$-algebra $A$, 
which is isometrically $A$-isomorphic to a nondegenerate left
 operator $A$-module, and suppose that $X$ is a dual 
Banach space.  Then the map $x \mapsto ax$ on $X$ is
$w^*$-continuous for all $a \in A$. 
\end{theorem}

\begin{proof}   Let $u$ be the map $x \mapsto ax$.
The adjoint of the canonical inclusion of the
predual of $X$ into its second dual, is  a
$w^*$-continuous 
contractive projection $q \colon X^{**} \to X$,
which induces an isometric  map
$v  \colon X^{**}/{\rm Ker}(q) \to X$.
By basic duality principles,
$v$ is $w^*$-continuous.
By the Krein-Smulian theorem (see Section 1),
it is a $w^*$-homeomorphism.  
We claim that: \begin{equation} \label{3mag}
q(u^{**}(\eta)) \; = \; u \, q(\eta) , \qquad \eta \in X^{**} .
\end{equation}
If  (\ref{3mag}) holds, then it follows that
$u^{**}$ induces a map $\dot{u}$ in
$B(X^{**}/{\rm Ker}(q))$, namely $\dot{u}(\dot{\eta}) =
(u^{**}(\eta))^{\dot{}}$, where $\dot{\eta}$ is the
equivalence class of $\eta \in X^{**}$ in the quotient.
Since $u^{**}$ is
$w^*$-continuous, so is $\dot{u}$,
by basic Banach space duality principles.
Using  (\ref{3mag}) it is easy to see
that $\dot{u} = v^{-1} u v$.
Since $\dot{u}, v,$ and $v^{-1}$ are
$w^*$-continuous, so is $u$, and thus the result is proved. 

In fact (\ref{3mag}) follows from Lemma \ref{newl}. 
Indeed, as in the proof of that result, $X^{**}$ may
be regarded as an operator $A^{**}$-module.  Therefore it is
an $A$-module, with module action $a \eta = u^{**}(\eta)$ in
the notation above, and $X$ is
an $A$-submodule. 
\end{proof}

\begin{definition}  A linear map $T$ on a Banach space will be
called a {\em Banach left multiplier} if there exists
a $C^*$-algebra $A$ (which may be taken
to be unital), a linear isometry $\sigma : X \to A$, and
an element $a \in A$, with $\sigma(Tx) = a \sigma(x)$ for all $x \in X$.
We say that that such a map $T$ is {\em adjointable}
 if also $a^* \sigma(X) \subset \sigma(X)$.
\end{definition}

\begin{corollary} \label{absm}  Every adjointable Banach left
multiplier of a dual
Banach space is $w^*$-continuous.
\end{corollary}

\begin{proof}
Let $T$ be as in the last definition, but with $X$ a dual Banach space.
Then $\sigma(X)$ is a left operator module over
the $C^*$-algebra generated by $1$ and $a$.
By the previous result,
 the map $L_a$ of left multiplication by $a$ is 
continuous in the
$w^*$-topology of $\sigma(X)$ induced by the predual of $X$.  It is then
clear that $T$ is $w^*$-continuous on $X$.
\end{proof}

As one application of the last 
results, we give  alternative
proofs of some important results about Hilbert $C^*$-modules
 from \cite{ACOT,OIAN}.
These particular proofs  will also be needed
 later in Section 5 (they generalize to the
nonselfadjoint situation).   
The reader who is only interested in 
our main results, and who is not familiar
with basic $C^*$-module theory, can skip the proof
of the next result, and also Section 5.
In fact we will not even give the basic definitions,
which may be found in any of the standard 
$C^*$-module texts, or \cite[Chapter 8]{BLM}.
We recall that  {\em ternary
rings of operators} (or {\em TROs}) are a
 simple example of  $C^*$-modules.
 They may be
taken to be spaces of the form $Z =
p A (1-p)$, for a unital $C^*$-algebra $A$ and a projection
$p \in A$ (e.g.\ see \cite{OIAN}
or (8.3) and the line above it in
\cite{BLM}).   For any TRO $Z$, $Z Z^*$ and $Z^* Z$ are
$C^*$-algebras, and $Z$ is a $Z Z^*$-$Z^* Z$-bimodule.
In fact every $C^*$-module may be represented 
canonically as a TRO (e.g. see 8.1.19 in \cite{BLM}).

\begin{corollary}  \label{eorz}  {\rm (Zettl, Effros-Ozawa-Ruan)} 
Let $Z$ full right Hilbert $C^*$-module over 
a $C^*$-algebra $B$, and suppose that 
$Z$ has a Banach space predual.  
If $N = M(B)$ then
$N$ and $B_N(Z)$ are $W^*$-algebras, the 
`inner product' on $Z$ is separately $w^*$-continuous
as a map into $N$,
 and $Z$ is a $w^*$-full
selfdual $W^*$-module over $N$ (see e.g.\ \cite[Section 8.5]{BLM} 
for definitions).  Moreover, $Z$ has a unique Banach space
predual, and this predual is also an operator space predual.
If $Z$ is a TRO with a Banach space predual, then
$Z$ is ternary isomorphic and $w^*$-homeomorphic
to a `corner' $q M (1-q)$, for  a $W^*$-algebra $M$ and a
projection $q \in M$.     
\end{corollary}        

\begin{proof}   We assume that 
$Z$ is a TRO, and $B = Z^* Z$.   This is purely for
notational simplicity, the  general case
is essentially identical.
The subalgebra $B_B(Z)$ of $B(Z)$ is $w^*$-closed
in the natural $w^*$-topology of $B(Z)$.
To see this, suppose that $(T_t)$ is a net in
$B_B(Z)$ converging in this topology to
$T \in B(Z)$.  Thus $T_t(y b) = T_t(y) b$
converges to $T(y b)$ in the $w^*$-topology of $Z$,
 for all $y \in Z, b \in B$.
On the other hand, $T_t(y)$ converges to $T(y)$.
Thus $T_t(y) b 
\to T(y) b$, by Theorem  \ref{opwe}.  It follows that
$T(y b) = T(y) b$, and so $T \in B_B(Z)$.
Thus $B_B(Z)$ is $w^*$-closed in $B(Z)$.

Let  $u \colon Z \rightarrow B$
be a bounded $B$-module map.
It is well-known (e.g.\ see 8.1.23 in \cite{BLM}) that we may choose
a contractive approximate identity $(e_t)_t$ for
$Z Z^*$,
with terms of the form $\sum_{k = 1}^n  x_k x_k^*$
for some $x_k \in Z$.   
Set $w_t = \sum_{k = 1}^n  x_k u(x_k)^*$ (which depends on $t$).
For $x \in Z$, 
\begin{equation} \label{7d}  u(e_t x) =
\sum_{k = 1}^n u(x_k) x_k^*  x  =
   w_t^*  x  . \end{equation}
 It follows   
that  $\Vert w_t \Vert^2 = \Vert u(e_t(w_t)) \Vert \leq 
\Vert u \Vert \Vert w_t \Vert$.  Thus
$(w_t)_t$ is a bounded net in
$Z$, and so it has a $w^*$-convergent subnet,
with limit $w$ say.  Replace the
net with the subnet.   By Theorem \ref{opwe},
 $z x^* w_t  \rightarrow z x^* w$, for all $x, z \in Z$.  Since
$u(e_t(x)) \rightarrow u(x)$  in norm, by
(\ref{7d}) we have $z u(x)^* =  z x^* w$,
 for all $x, z \in Z$.   Thus $u(x) = w^* x$,
and so $Z$ is selfdual over $B$.
It follows immediately that $B_B(Z)$ is the $C^*$-algebra of
`adjointable' maps on $Z$ (e.g.\ see 8.5.1 (2) in \cite{BLM}).
Equivalently, by a result of Kasparov 
(e.g.\ see 8.1.16 in \cite{BLM}),  $B_B(Z) = M(Z Z^*)$.  Since $B_B(Z)$ has a 
predual, it is a $W^*$-algebra in its
natural $w^*$-topology (that is, a bounded net
of maps converges if and only if they converge
as maps on $Z$, in the point-$w^*$-topology).    By symmetry,
$N = M(Z^* Z)$ is a  $W^*$-algebra in a topology
for which a bounded net $(n_t)$ converges to $n$
if and only if $z n_t \to z n$ for all $z \in Z$.

Claim: the inner product 
is separately $w^*$-continuous.  
Suppose that $(y_t)$ is a bounded net 
in $Z$ converging in the $w^*$-topology of $Z$ to
$y \in Z$, and that $w \in Z$ is fixed.
Suppose that $(w^* y_{t_\mu})$ is a subnet of 
$(w^* y_{t})$, with $w^*$-limit $n \in N$.
If $z \in Z$, then $(z w^* y_{t_\mu})$ converges
both to $z w^* y$ and to $z n$, by Theorem \ref{opwe} and the 
fact at the end of the last paragraph.  It follows that 
$n = w^* y$.  Thus the claim is proved (using also 
the Krein-Smulian theorem as mentioned in Section 1). 

The other assertions of the Theorem now all follow immediately from          
standard facts about selfdual modules (e.g.\ see 8.5.1-8.5.4, and 
8.5.10, in \cite[Section 8.5]{BLM}).  These facts are also all
mentioned in Section 5 below, in a more general setting
(see particularly Lemma \ref{7prez}).
\end{proof}
          
\section{Multipliers and duality}

\begin{theorem} \label{3lawc} 
Every left multiplier of a dual operator space
is $w^*$-continuous.  \end{theorem}
 
\begin{proof}   If $u \in \M_l(X)$, then
$u^{**} \in  B(X^{**})$.  As in Theorem \ref{opwe},
let $q \colon X^{**} \to X$
be the canonical projection, which now
is completely
contractive.  As in that result, it suffices to show
that
 \begin{equation} \label{3mageq}
q(u^{**}(\eta)) \; = \; u \, q(\eta) , \qquad \eta \in X^{**} .
\end{equation}
In order to prove 
(\ref{3mageq}), we let $Z$ be an injective
envelope of $X$, viewed
 as a TRO $p D (1-p)$, 
for a unital $C^*$-algebra $D$ and a projection $p \in D$. 
(see e.g.\ \cite[Sections 4.2 and 4.4]{BLM} or
\cite[Chapter 16]{CBMAO}). 
If $E = Z^{**}$ then $E  = p D^{**} (1-p)$ is also a TRO.
Clearly $X^{**}$ may be regarded as
a $w^*$-closed subspace of $E$, and thus by injectivity 
of $Z$, we can extend the map $q$   above, to a 
completely contractive map $\theta : E \to Z$.
Since $\theta_{\vert X} = I_X$, by the rigidity property
of the  injective
envelope we must have $\theta_{\vert Z} = I_Z$.
Thus $\theta$ is a completely contractive 
projection from $E$ onto 
$Z$.  By Lemma \ref{newl}, $\theta$ is a left $p D p$-module map.
Let  $a \in pDp$ be such that  $a x = u x$ for all
$x \in X$, as in \cite[Theorem 4.5.2]{BLM} or
\cite[Chapter 16]{CBMAO}.    
Since $\theta$ is a left $pDp$-module map, 
 \begin{equation} \label{eq1}
\theta(a \eta) = a \theta(\eta) = a q(\eta) , \qquad \eta \in X^{**}.
\end{equation} On the other hand, we claim that
\begin{equation} \label{eq2}
a \, \eta \; = \; u^{**}(\eta) , \qquad
\eta \in X^{**} . \end{equation} 
To see this, view both sides as functions from $X^{**}$ into $E$.
Then both functions are $w^*$-continuous
(note that since $E = p D^{**} (1-p)$, left 
multiplication by the element $a \in p D p \subset p D^{**} p$ is
$w^*$-continuous).
On the other hand, (\ref{eq2})
 certainly holds if $\eta \in X$,
and by $w^*$-density it must therefore hold for $\eta \in X^{**}$.
By (\ref{eq2}), we have that $\theta(a \eta) = \theta(u^{**}(\eta))
= q(u^{**}(\eta))$.  This together with
(\ref{eq1}) proves that
$q(u^{**}(\eta))=aq(\eta)=uq(\eta)$, which is
(\ref{3mageq}).  
\end{proof}

The last theorem has very many applications:

\begin{corollary} If $X$ is a dual Banach space, then every
multiplier of $X$ in the sense of the second paragraph
of our paper, 
is $w^*$-continuous.
\end{corollary}
 
\begin{proof}   This follows from Theorem \ref{3lawc} applied to
${\rm Min}(X)$, which is a dual operator space by 1.4.12
in \cite{BLM}.
By 4.5.10 in \cite{BLM},  $\M_l({\rm Min}(X))$ is the set of 
multipliers in the sense of \cite[Section I.3]{MIBS}.
\end{proof}
      
The last result was known for `centralizers' of Banach spaces
(e.g.\ see \cite[Theorem 1.3.14]{MIBS}),
but seems not to be known for the larger class of Banach space
multipliers.   We thank E. Behrends for confirming this.
We imagine that it may be useful in that theory too.   

The following answers a question that has also been open for many
years:
 
\begin{corollary} \label{3lawc2}  If $B$ is an operator algebra
which is also a dual operator space, then the product on $B$
is separately $w^*$-continuous.
\end{corollary}
 
\begin{proof}  If $a\in B$,
 then the map $b \mapsto ab$ on $A$ is clearly a left multiplier,
and therefore is $w^*$-continuous by Theorem \ref{3lawc}.  Similarly
the product is $w^*$-continuous in the first variable.  \end{proof}

Putting Corollary \ref{3lawc2} together with the main result in 
\cite{AOSC}, we obtain the following improved characterization
of $\sigma$-weakly closed operator algebras:
 
\begin{corollary} \label{3lawc3}  If $B$ is an operator algebra
which is also a dual operator space, then $B$ is completely 
isometrically isomorphic, via a $w^*$-homeomorphic homomorphism,
to a $\sigma$-weakly closed subalgebra of $B(H)$, for 
some Hilbert space $H$.
Conversely, every $\sigma$-weakly closed subalgebra of $B(H)$
is a dual operator space.
\end{corollary} 

Corollaries \ref{3lawc2} and 
 \ref{3lawc3} were obtained in \cite{MADO} in the case that 
$B$ also has an identity of norm $1$. 
Corollary \ref{3lawc3} makes the following definition
quite substantive:
 
\begin{definition} \label{3lawc4}  A {\em dual operator algebra}
is an operator algebra which is also a dual operator space.
\end{definition}

The following was noticed together with Le Merdy:

\begin{corollary} 
The Arens product is the only operator algebra  product
on the second dual of an  operator algebra $A$, which extends 
the product on $A$. 
\end{corollary}  

\begin{proof}  By Corollary \ref{3lawc2}, any such 
product is separately $w^*$-continuous.
\end{proof} 

\begin{corollary} Every quasimultiplier (in the sense of
\cite{KP}) of a dual operator space, is separately $w^*$-continuous.
\end{corollary}

\begin{proof}   This follows from Corollary \ref{3lawc2}, and
Kaneda's correspondence 
between contractive quasimultipliers and operator algebra
 products.  \end{proof}

\begin{corollary}  Suppose that $B$ is a operator algebra
with a bounded approximate identity, and with an
operator space predual.  Then $B$ has an
identity (of norm possibly $> 1$).
\end{corollary}
 
\begin{proof}   If $e$ is a $w^*$-limit of a bounded 
approximate identity, then $e$ is an
identity by Corollary \ref{3lawc2}.   \end{proof}
                           
We do not know if the last several results are true
if we replace `dual operator space' by `dual Banach space',
and we shall establish below 
some partial results along these lines.
  Corollary \ref{3lawc3}  is not true, as it stands, 
in this case (see e.g.\ \cite{AOSC,MADO} for counterexamples).
  However if Corollary \ref{3lawc2} were true in this case, 
then the proof of 
the main result of \cite{AOSC} would yield that
every  operator algebra with a Banach space predual,
is {\em isometrically}
 isomorphic, via a $w^*$-homeomorphic homomorphism,
to a $\sigma$-weakly closed subalgebra of $B(H)$.  

We recall that $\M_l(X)$ contains the $C^*$-algebra $\A_l(X)$
of {\em left adjointable multipliers}.  These are 
the left multipliers as defined in the second paragraph
of our paper in terms of a complete isometry $\sigma$,
 but also satisfying $a^* \sigma(X) \subset \sigma(X)$
in the language of that definition.  It is shown in 
\cite{MADO,OMAM} that if $X$ is a dual operator space
then $\M_l(X)$ is a dual operator algebra,
$\A_l(X)$ is a $W^*$-algebra, and every
$T \in \A_l(X)$ is $w^*$-continuous.  These results
have been key to  work on multipliers and
$M$-ideals following \cite{OMAM} (see e.g.\
\cite{MCAA,TCOO}).   We shall see that
one of these results is true, and the others false,
if $X$ merely has a Banach space predual.   Indeed 
from Corollary \ref{absm} we have:  

\begin{corollary}  \label{adjws}  Let $X$ be an
operator space, which is a dual Banach space.
Then every $T \in \A_l(X)$ is $w^*$-continuous.
\end{corollary}

\begin{proposition} \label{delt} Let $B$ be an approximately unital
operator algebra which is
a dual Banach space.  Then $B$ is unital,
$\Delta(B) = B \cap B^*$ is a $W^*$-algebra,
and if $b \in \Delta(B)$ then the maps $a \mapsto ab$ and
$a \mapsto ba$ are $w^*$-continuous.
\end{proposition} 
 
\begin{proof}    It is shown in 
\cite[Theorem 2.5]{MADO} that $B$ is
unital.  Let $A = B^{**}$ and let $q : A \to B$ 
be the canonical projection.   Since $q(1) = 1$, 
$q$ takes Hermitian elements to Hermitian elements.
That is, $q$ induces a $w^*$-continuous projection of 
$\Delta(A)$ onto $\Delta(B)$.    Thus 
$\Delta(B)$ is isometric to the dual space 
$\Delta(A)/{\rm Ker}(Q)$, and so 
$\Delta(B)$ is a $W^*$-algebra.  
The last part follows from e.g.\ Corollary \ref{absm},
but we give an independent proof:
 Any projection
in $\Delta(B)$ corresponds to a left $M$-projection
 on $B$.  The latter is necessarily 
$w^*$-continuous, by Proposition 3.3 in
\cite{OMAM}.  Since the projections
densely span a $W^*$-algebra, left multiplication 
by any element in $\Delta(B)$ is $w^*$-continuous.
A similar argument pertains to right multiplications.  
  \end{proof} 

The following is in stark contrast to the `dual operator space
case' mentioned above Corollary \ref{adjws}.  It also shows 
that at least one plausible variant of the second assertion
in Proposition \ref{delt} fails
for general operator spaces:

\begin{proposition} \label{naw}  There exists an operator system
$X$ which is
a dual Banach space, but for which $\A_l(X)$ is not a 
$W^*$-algebra (or even an $AW^*$-algebra), and
$\M_l(X)$ is not a dual Banach space.
\end{proposition}

\begin{proof}   Let $X$ be the operator system in Proposition
\ref{nonp}.
   We will show that $\A_l(X) = \M_l(X) \cong
S^\infty + \Cdb I$.   To this end,
we first claim that $D = B(H) \oplus^\infty Q^{\rm op}$ 
is the $C^*$-envelope of $X$.   Let $B$ be the 
$C^*$-algebra generated by $X$ in $D$.
  If $(a,\dot{a}), (b,\dot{b})
\in X$, then   $(a,\dot{a}) (b,\dot{b}) =
(ab,\dot{a} \dot{b}) \in B$.  Since $(ba, \dot{a} \dot{b}) \in X$,
we have $(ab - ba,0) \in B$.   If $(c,\dot{c}), (d,\dot{d}) \in X$,
then $(c,\dot{c}) (ab - ba,0) (d,\dot{d}) = (c(ab - ba)d,0) \in B$.
However $\{ c(ab - ba)d : a, b, c, d \in B(H) \}$ densely spans $B(H)$,
and so $B$ contains $B(H) \oplus 0$.  Since also $X \subset B$,
we have $0 \oplus Q^{\rm op}  \subset B$, and  
it follows that $B = D$.
Thus $X$ generates $D$ as a
$C^*$-algebra.   To see that $D$ is
the $C^*$-envelope of $X$, suppose that $J$ were a nontrivial ideal in
$D$ such that the canonical map $D \to D/J$ is 
completely isometric on $X$.  Since $B(\ell^2)$ has only one nontrivial
closed ideal, and therefore $Q$ has none,
 $J$ must be one of the four spaces  $B(H) \oplus 0, S^\infty
\oplus 0,  0 \oplus Q^{\rm op}, S^\infty
\oplus Q^{\rm op}.$  Thus
$D/J$ is of the form $0 \oplus Q^{\rm op},
Q \oplus Q^{\rm op}, B(H) \oplus 0,
$ or  $Q \oplus 0$.
In any of these cases we obtain a contradiction.
For example, the third case yields a contradiction, because,
 by the discussion in the proof of 
Proposition
\ref{nonp}, the map 
$(a,\dot{a}) \to (a,0)$ from $X$ to $B(H) \oplus 0$
is not a complete
isometry.  This proves the claim.     

As explained in Section 1,
$\M_l(X) = \{ a \in D : a X \subset X \}$.  Clearly
$(b + \lambda 1,\lambda \dot{1}) \in \A_l(X)$, for 
any $b \in S^\infty$.
Since $1 \in X$, if $a \in \M_l(X)$ then $a = (b,\dot{b})$ for a
$b \in B(H)$ such that $(bc,\dot{b} \dot{c}) \in X$ for all
$c \in B(H)$.  That is, $\dot{b} \dot{c} =  \dot{c} \dot{b}$,
so that $\dot{b}$ is in the center of $Q$.  However 
the center of $Q$ is trivial (see e.g. \cite{JP}, we thank V. Zarikian
for communicating this reference to us).  Thus $b \in 
S^\infty + \Cdb 1$.  Hence $\M_l(X) = 
\A_l(X) \cong S^\infty + \Cdb 1$.  
 \end{proof}  

Theorem \ref{3lawc} will also be an important tool
for future work on  {\em operator modules}.
For example, in \cite{MADO} we were able to 
improve in several ways on a theorem of Effros and 
Ruan characterizing
certain operator modules over von Neumann
algebras \cite{ROOB}.    Theorem \ref{3lawc} allows precisely the same
improvements for `normal dual operator modules' over 
unital dual operator algebras.  
In particular, Theorem \ref{3lawc} 
 shows that the {\em left normal} hypothesis
used in \cite{MADO} is automatic,
and may therefore be removed.
We state a sample of other consequences:

\begin{corollary}  \label{lawcmo} Suppose that $X$ is a 
left operator $A$-module, where $A$ is approximately 
unital, and suppose that $X$ is
also a dual operator space.  Then for any $a \in A$,
the map $x \mapsto ax$ is automatically $w^*$-continuous.
This is also true, if $X$ merely is a
dual Banach space, providing that $A$ is
a $C^*$-algebra.
  \end{corollary}

This corollary allows one to eliminate one of the 
hypotheses in the well-known definition of a
{\em normal dual operator bimodule} (e.g.\ see
\cite{ROOB}).   Thus we may define, for example,
a {\em left normal dual operator module}
to be a left operator module $X$ over a dual operator algebra $M$,
such that $X$ is a dual operator space, and
the module action $M \times X \to X$ is $w^*$-continuous in the first variable.
           Similar definitions hold for right modules and bimodules.
                
\begin{corollary} \label{xisnd}
Let $X$ be a dual operator space.  Then $X$ is a 
normal dual $\M_l(X)$-$\M_r(X)$-bimodule.
\end{corollary}
  
Conversely, any normal dual operator module or bimodule action on a
dual operator space $X$ `factors  through' the
one in Corollary \ref{xisnd}; and moreover there is a tidy 
`representation theorem' for such modules.
For details, see \cite{MADO} or 4.7.6 and 4.7.7 in \cite{BLM}.

\bigskip

The following is a Banach module characterization of $w^*$-closed subspaces
$X$ of $B(K,H)$ which are invariant under the action of two 
$W^*$-algebras $M$ and $N$ on $H$ and $K$ respectively
(that is, $\pi(M) X \subset X$ and $X \theta(N)  \subset
X$, where $\pi$ and $\theta$ are normal $*$-representations of
  $M$ and $N$).   Our theorem is
the Banach module variant of an earlier operator module 
characterization of such bimodules due to Effros and Ruan \cite{ROOB};
however it has potential to be even more useful in certain contexts
since our condition (ii) is easier to verify.
  Before we state the theorem we give a definition: to say that
the unit ball $Ball(X)$ is $M$-$N$-{\em absolutely convex} is to say that
$\Vert \sum_{k=1}^m  m_k x_k n_k \Vert \leq 1$ whenever 
$x_1, \cdots, x_m \in Ball(X)$, and $m_1, \cdots, m_m \in M, 
n_1, \cdots, n_m \in N$ with
$\Vert \sum_{k=1}^m  m_k m_k^* \Vert \leq 1$ and  $\Vert \sum_{k=1}^m  
n_k^* n_k \Vert \leq 1$.

\begin{theorem}  \label{ERbm}  Let $M$ and $N$ be 
$W^*$-algebras, and let $X$ be a Banach $M$-$N$-bimodule
(we assume that $1_M x = x 1_N = x$ for all $x \in X$).
Suppose that $X$ is also a dual Banach space.  
The following are equivalent:
\begin{itemize} \item [{\rm (i)}]  
There exist Hilbert spaces $K$ and $H$, 
 a $w^*$-continuous isometry $\Phi : X \to B(K,H)$,
and normal $*$-representations $\pi$ and $\theta$ of $M$ and $N$
on $H$ and $K$ respectively, such that
$\Phi(m x n) = \pi(m) \Phi(x) \theta(n)$ for $x \in X, m \in M, n \in N$;
\item [{\rm (ii)}]  The unit ball of $X$ is $M$-$N$-absolutely convex,
and for all $x \in X$ the canonical maps $M \to X$ and $N \to X$ given
by $m \mapsto m x$ and $n \mapsto x n$, are $w^*$-continuous;
 \item [{\rm (iii)}]  The unit ball of $X$ is $M$-$N$-absolutely convex,
and the bimodule action $M \times X \times N \to X$ is
separately $w^*$-continuous.
\end{itemize}
\end{theorem}

\begin{proof}  (iii) $\Rightarrow$ (ii) \ This is trivial.

(ii) $\Rightarrow$ (iii) \  The condition implies by e.g.\
\cite[Theorem 2.1]{TMOM} that there is an operator space structure
on $X$ for which $X$ becomes an operator $M$-$N$-bimodule.
Now (iii) is clear from Corollary \ref{absm}.

(iii) $\Rightarrow$ (i) \ As in the lines above,
$X$ may be viewed as an operator $M$-$N$-bimodule.  By 
e.g.\ 3.8.9 in \cite{BLM}, $X^{**}$ is
an operator $M$-$N$-bimodule too.   As in the first
few lines of 
Theorem \ref{opwe},  there is a 
canonical $w^*$-continuous contractive projection $q \colon X^{**} \to X$,
which induces an isometric $w^*$-homeomorphism  
$v  \colon X^{**}/{\rm Ker}(q) \to X$.   By 
Lemma \ref{newl}, $q$ is an $M$-$N$-bimodule map,
and therefore so also is $v$.  
Indeed, $X^{**}/{\rm Ker}(q)$ is an operator $M$-$N$-bimodule
isometrically $M$-$N$-isomorphic to $X$, via $v$.
We assign $X$ a new operator space structure so that
$v$ becomes a complete isometry.  Since 
$X^{**}/{\rm Ker}(q)$ has an
operator space predual (namely ${\rm Ker}(q)_\perp$),
so now does $X$.  Moreover, we have not changed the
$w^*$-topology on $X$, since $v$ was a $w^*$-homeomorphism
originally.  Now $X$ is a normal dual 
 operator $M$-$N$-bimodule, and hence 
we obtain the desired representation from e.g.
\cite{ROOB} or \cite[Theorem 3.8.3]{BLM}.     

 (i) $\Rightarrow$ (ii) \  This is the routine `easy direction' of
such theorems, and is left here as an exercise.  
\end{proof}

The next section will continue to demonstrate that
operator space multipliers, and Theorem \ref{3lawc},
will be key to future studies of operator modules.

\section{Nonselfadjoint generalization of $W^*$-modules}

Notions of Morita equivalence appropriate to
nonselfadjoint operator algebras, and of `rigged modules',
were developed in the last ten years 
in \cite{COOM,AGOH}.  These notions generalize the `strong
Morita equivalence' of $C^*$-algebras due to Rieffel,
and the `$C^*$-modules' used heavily in that theory.
There is a parallel theory, mainly due to Paschke
and Rieffel (see 
\cite{IPMO,MEFC} or \cite[Section 8.5]{BLM})
appropriate to $W^*$-algebras: the corresponding
notions are
sometimes called `$W^*$-algebra Morita equivalence',
and  `$W^*$-modules'. 
  Hitherto there has been no attempt in the 
literature to generalize this `weak' version of the 
theory, to nonselfadjoint 
dual operator algebras.  One main reason 
for this, we believe, is that the technical tools
were not all available or fully developed.    It seems that 
operator space multipliers and Theorem \ref{3lawc},
were one of the missing ingredients in getting this
theory started.  We can show 
that with the addition of this ingredient,
one can obtain a theory that 
generalizes several important aspects of the $W^*$-algebra case.
At the same time, this will illustrate how
Theorem \ref{3lawc} may be powerfully used in practice.
Our intention is to be very brief; the reader will need to 
consult the papers \cite{COOM,AGOH} for additional
definitions and details.  
 
In the following discussion,
 $Y$ is a right $M$-rigged module, in the sense 
of  \cite{AGOH}, over an approximately unital operator algebra
$M$.   Then there is a canonical left $M$-rigged module
$X = \tilde{Y}$,
and a canonical pairing $(\cdot,\cdot) : X \times Y \to M$
 (see \cite{AGOH} or \cite[Chapter 4]{COOM}).  In 
our case, $M$ will usually be a dual operator algebra.
  We say that 
$Y$ is {\em selfdual} over $M$, if every completely
bounded $M$-module map $f : Y \to M$ is of the form
$(x,\cdot)$ for a fixed $x \in X$, and every
completely bounded
$M$-module map $g : X \to M$  is of the form
$(\cdot,y)$ for a fixed $y \in Y$.  
If $Y$ is  selfdual then every completely 
bounded $M$-module map from $Y$ into another rigged
 $M$-module $Z$
is {\em adjointable}.  This follows by considering
the $M$-valued $M$-module map $(w,u(\cdot))$ on $Y$,
for fixed $w \in \tilde{Z}$, just as in the
$C^*$-module case (e.g.\ see 8.5.1 (2) in \cite{BLM}).
Indeed, the proofs of the next two results are also
essentially just as in Section 8.5 of
\cite{BLM}, simply replacing appeals to 
$C^*$-module facts by appeals to the matching results
for rigged modules from \cite{AGOH,COOM}.  Thus
we omit essentially all of these proofs. 
 
\begin{lemma} \label{7prez} Let $Y$ be a 
right rigged $M$-module over a 
unital dual operator algebra $M$.  Then:
\begin{enumerate}
\item [{\rm (1)}]
$Y$  is a selfdual rigged $M$-module if and only if
$X$ and $Y$ have  Banach space preduals
 with respect to which $(\cdot,\cdot)$  
is separately $w^*$-continuous.
\end{enumerate}
If $Y$ is a  selfdual rigged $M$-module, then:
\begin{enumerate}
\item [{\rm (2)}]
$X$ and $Y$ have unique Banach space preduals
with respect to which $(\cdot,\cdot)$  
is separately $w^*$-continuous.
\item [{\rm (3)}]  With respect to the $w^*$-topology induced
by the predual in {\rm (2)},
a bounded net
$(y_t)_t$ converges to $y$ in $Y$
if and only if
$(x , y_t ) \rightarrow (x , y)$
in the $w^*$-topology of $M$, for all $x \in X$.
Similarly for  bounded nets in $X$.
  \item [{\rm (4)}]
Let $W =  M_* \overset{\frown}{\otimes}_M X$ and 
$Z = Y \overset{\frown}{\otimes}_M  M_*$ (see \cite[Section 3.4]{BLM}).
  Then $W$
and $Z$ are  operator space preduals of $Y$  and $X$ 
respectively, inducing the
$w^*$-topology in {\rm (2)} and {\rm (3)} above.  
\item [{\rm (5)}]  
The canonical map $m \mapsto y m$ from 
$M$ to $Y$ is $w^*$-continuous in the topology in 
{\rm (3)}, for all fixed $y \in Y$.
  \end{enumerate}
\end{lemma}

\begin{proof}  We will simply prove (5), which was not
mentioned in the matching result from \cite{BLM}.
If $(m_t)$ is a bounded net 
converging to $m$ in the $w^*$-topology of $M$, and if $x \in X, y \in Y$,
then  we have $(x , y m_t) = (x , y) m_t \to (x , y) m
= (x , y m)$, by the separate $w^*$-continuity 
of the product in $M$.  Thus by (3),
$y m_t \to y m$.  The result follows by the 
Krein-Smulian theorem (see Section 1).   \end{proof}
 
We will henceforth use the phrase
{\em the $w^*$-topology} of a selfdual $M$-rigged module, for
the (unique) topology in (2)--(4) above.

\begin{corollary} \label{7prez2}
Suppose that $Y$ is a selfdual right $M$-rigged module
over a unital dual operator algebra
$M$.  Then: \begin{enumerate}
\item [{\rm (1)}]    $CB_M(Y) = \Bdb_M(Y)$, the operator
algebra of `adjointable' $M$-module maps, and
this is a dual operator algebra.
\item [{\rm (2)}]  A bounded
net $(T_i)_i$ in $CB_M(Y)$ converges in the $w^*$-topology
 to $T \in B_M(Y)$
if and only if $T_i(y) \to T(y)$ in the
$w^*$-topology of $Y$, for all $y  \in Y$.
Indeed, $Y \overset{\frown}{\otimes}_M W$ is
a predual for $CB_M(Y)$,
where $W$ is as in Lemma {\rm \ref{7prez} (4)}.
\end{enumerate}
\end{corollary}

Similarly it follows, as in \cite[Corollary 8.5.8]{BLM},
that any bounded $M$-module map between
selfdual right rigged $M$-modules, is $w^*$-continuous. 
 
For a right rigged module $Y$ over an operator
algebra 
$A$, we will consistently write
${\mathcal I}$ for the closed span of the range
of the canonical pairing $(\cdot,\cdot)$ in $A$.
We say that $Y$ is {\em full} over $A$, if $A = 
{\mathcal I}$.
If $A$ is a dual operator algebra, we write
${\mathcal I}^w$ for the $w^*$-closure of this
span, and say that $Y$ is
{\em $w^*$-full} if ${\mathcal I}^w = A$.  
In general though, ${\mathcal I}$ and 
${\mathcal I}^w$ are both ideals in $A$.             
Henceforth, we say that a right rigged module 
$Y$ is {\em a (right) rigged-equivalence module}, if 
${\mathcal I}$ has a contractive approximate identity,
and the canonical map $X \otimes_{h} Y \to 
{\mathcal I}$ is a complete
quotient map.  This is equivalent to saying that 
${\mathcal I}$   possesses a contractive approximate 
identity of a certain special form, or to saying that
$Y$ is a strong Morita equivalence
$\Kdb_A(Y)$-${\mathcal I}$-bimodule.  E.g.\ see \cite{COOM,AGOH}
for more details.  In this case, and if also $A$
is a dual operator algebra, then
by considering a $w^*$-limit
of the contractive approximate
identity, it follows that ${\mathcal I}^w$
is unital.  

The property of
selfduality defined earlier
does not depend essentially on $M$: that
is, $Y$ is selfdual over $M$ if and only if
$Y$ is selfdual over ${\mathcal I}$ or over 
the multiplier algebra $M({\mathcal I})$.  The proof of this
is identical to  \cite[Lemma 8.5.2]{BLM}. 

We will write $LM(A)$ and $RM(A)$ for the left
and right  multiplier
algebras of $A$.  For example, $LM(A)$ may be identified with
$CB_A(A)$ (see e.g. \cite[Section 2.6]{BLM}).

\begin{lemma} \label{wsid}  Let $J$ be a $w^*$-dense
norm-closed two-sided 
ideal in a dual operator algebra $M$, and suppose that
$J$ is approximately unital.  Then $M$ is the multiplier
algebra $M(J)$, and the latter equals $LM(J)$ and $RM(J)$.
\end{lemma}

\begin{proof}  (Sketch) \  In fact this works more generally in 
the setting of Banach algebras, provided that the
product on $M$ is separately $w^*$-continuous.
There is a canonical complete contractive homomorphism
$M \to CB_J(J)$, and the latter space
is just $LM(J)$.   This map is 1-1 by the $w^*$-density
of $J$, that it is completely isometric
and surjective is easily seen by considering,
for any $T \in CB_J(J)$, a $w^*$-limit point of $(T(e_t))$ in
$M$, where $(e_t)$ is the approximate identity for $J$.  
The other assertions are now easy.       
\end{proof}             
 
\begin{corollary} \label{7jism}
Let $Y$ be a rigged-equivalence module,
over a dual unital operator algebra $M$.  In the notation above,
${\mathcal I}^w$ is the multiplier algebra
of ${\mathcal I}$. 
\end{corollary}
 
\begin{proof}  Clearly  ${\mathcal I}$ is a $w^*$-dense
ideal in ${\mathcal I}^w$.   
\end{proof}

We now seek to generalize Zettl's theorem 
(cf.\ Corollary \ref{eorz}) to rigged modules.  
It is natural to assume in our context that 
$Y$ and $X = \tilde{Y}$ both have
an operator space predual.   Our main theorem,
Theorem \ref{3lawc}, then yields the following
corollary, which in turn will yield the nonselfadjoint analogue of Zettl's theorem.
 
\begin{corollary} \label{fromm}  Let $Y$ be a right rigged
module over an approximately unital 
operator algebra $A$, and suppose that $Y$ and $X = \tilde{Y}$ both have
an operator space predual.  Then the maps $y \mapsto y' (x',y)$,  
$x \mapsto (x,y') x'$, $y \mapsto y a$,
and $x \mapsto a x$, are automatically 
$w^*$-continuous on $Y$ and $X$ respectively,
for all fixed $x',x'' \in X$, 
$y',y'' \in Y$, and $a \in A$. 
\end{corollary}

\begin{proof}   From the
theory of rigged modules, it is clear that these
maps are operator space multipliers.  For example, the
first of  these maps belongs to $B = \Kdb_A(Y)$,
and $Y$ is an operator $B$-$A$-bimodule.
Thus we can appeal to 
Corollary \ref{lawcmo} to see that this map,
and also the map $y \mapsto ya$,
are $w^*$-continuous on $Y$.  Similarly for
$X$.   \end{proof}

\begin{theorem}  \label{7osz}  Let $Y$ be a full
right rigged-equivalence module over an
approximately unital operator algebra $A$, 
and suppose that $Y$ and $X = \tilde{Y}$  are
dual operator spaces.   If $M = M(A)$ then
$M$ and $\Bdb_A(Y)$ are dual operator 
algebras, and $Y$ is a $w^*$-full
selfdual $M$-rigged module.
\end{theorem}
 
\begin{proof}  We follow the proof of 
Corollary \ref{eorz}.  As in that proof,
 but also using Corollary \ref{fromm}, 
$CB_A(Y)$ is a $w^*$-closed
subalgebra of $CB(Y)$.  
From the theory of strong Morita equivalence
(see e.g.\ \cite[Theorem 4.9]{COOM}),
 $CB_A(Y)$ is an operator algebra,
hence it is a dual operator algebra, by Corollary \ref{3lawc3}.
Similarly for $_ACB(X)$.
  
If $B = \Kdb_A(Y)$, then from the theory of strong Morita equivalence 
$LM(A) \cong CB_B(X)$, and $RM(A) \cong \, _BCB(Y)$, completely 
isometrically.  This may be seen from the fact that 
strong Morita equivalence implements a `completely isometric'
equivalence between the categories of right modules over 
$A$ and $B$ (see \cite[p. 25]{COOM}), thus
$$LM(A) = \, CB_A(A) \cong \, CB_B(A \otimes_{hA} X) \cong \,
 CB_B(X) . $$   The map here from $LM(A)$ into 
$CB(X)$ may be checked to be the canonical one: if $\eta \in LM(A),
a \in A, x \in X$, then $\eta$ takes $a x$ to $(\eta a) x$. 
Similarly, for the map from $RM(A)$ into $CB(Y)$.
Thus we identify the operator algebras
$LM(A)$ and $RM(A)$ with 
$CB_B(X)$ and $_BCB(Y)$ respectively; and by the
 argument above, these subspaces are dual operator algebras,
and $w^*$-closed subspaces of $CB(X)$ and $CB(Y)$ respectively.
Let $u \in CB_A(Y,A)$.
  Following  Corollary \ref{eorz},
we choose a contractive approximate
identity $(e_t)_t$ for $\Kdb_M(Y)$,
of the form $\sum_{k = 1}^n [ y_k , x_k]$ for some $x_k \in X,
y_k \in Y$ as in
e.g.  \cite[Theorem 5.2]{AGOH}.   For $y \in Y$, we have
\begin{equation} \label{7defwt}  u(e_t(y)) =
\sum_{k = 1}^n u(y_k) ( x_k , y ) =
( \sum_{k = 1}^n  u(y_k) x_k , y )
=  ( w_t , y )  ,  \end{equation}
where  $w_t = \sum_{k = 1}^n  u(y_k) x_k$.
Using the fact that the  canonical map
$X \to CB(Y,A)$ is an isometry (see e.g.\ \cite[Theorem 4.1]{COOM}), 
we have $\Vert w_t \Vert = \Vert u \circ e_t \Vert_{cb} \leq \Vert u
 \Vert_{cb}$.  Thus
 $(w_t)_t$ is bounded in
$X$, and we can proceed as in Corollary \ref{eorz},
but also using Corollary \ref{fromm},
to find $w \in X$ with  
$u(y) x' = (w , y ) x'$ for all $x' \in X, y \in Y$, so that 
$u(y) = (w , y )$.  Similarly,
any 
$A$-valued $A$-module map on $X$,
is given by $(\cdot,y)$ for a fixed $y \in Y$. 
Thus $Y$ is selfdual as an $A$-module.   
It follows, as asserted earlier,
that $CB_A(Y) = \Bdb_A(Y)$.   
As in the first centered equation of the present proof,
we have $CB_A(Y) \cong LM(\Kdb_A(Y))$.
This isomorphism carries $\Bdb_A(Y)$ onto $M(\Kdb_A(Y))$,
as one may check somewhat analogously to the proof
of 8.1.16 in \cite{BLM} (see \cite[Theorem 3.8]{AGOH}). 
We have now shown that $M(\Kdb_A(Y)) = LM(\Kdb_A(Y))$,
and it follows by symmetry that $M(A) = RM(A)$.
Similar arguments involving $X$ show that $M(A) = LM(A)$.

Now $Y$ is selfdual over $M(A)$ too,
as remarked above Lemma \ref{wsid}.
By Corollary \ref{7jism}, it is clear that $Y$ is a 
$w^*$-full module over $M(A)$.       \end{proof}

\begin{corollary}  \label{7masd}
Let $Y$ be a right rigged-equivalence
module over a dual operator algebra $M$.
Then $Y$  is a selfdual
 $M$-rigged module if and only if $Y$ and 
$X = \tilde{Y}$ possess operator space preduals.
\end{corollary}

\begin{proof}  By Lemma \ref{7prez} we need
only prove one direction.  Assuming
the existence of operator space preduals,
by the previous result and Lemma \ref{7jism}, $Y$ is
selfdual over  ${\mathcal I}^w
= M({\mathcal I})$.  It follows as in \cite[Lemma 8.5.2]{BLM}
that $Y$ is selfdual over  $M$.
\end{proof}

{\bf Remark.}  The wary reader may wonder whether
the given preduals in Corollary \ref{7masd} induce the
$w^*$-topology mentioned after 
Lemma \ref{7prez}.  Unlike the 
$W^*$-algebra case, in fact they may not, if these 
preduals were chosen poorly.  This is clear
by considering the simplest example: $A = B = X = Y$, a
unital  operator algebra
with several unrelated preduals (e.g.\ see 
\cite{MADO} or \cite[Corollary 2.7.8]{BLM}). 
 There are other conditions
one may impose,  that will alleviate this situation.
For example, if one also insists in Corollary \ref{7masd}
that $X$ and $Y$ be normal dual $A$-modules (defined
above Corollary \ref{xisnd}).  We leave this as an exercise for the interested reader
(see the ideas in the 
proof of Proposition \ref{newp} below).

\medskip 
 
Let $M$ and $N$ be 
two unital dual operator  algebras,
and suppose that there exist
$w^*$-dense norm closed ideals of $M$ and $N$ respectively,
which are strongly Morita equivalent in the
sense of \cite{COOM}, via 
equivalence bimodules $X$ and $Y$.  
By Lemma \ref{wsid}, $M$ and $N$ are
the multiplier algebras of these ideals.
Hence $Y$ is canonically
an  operator $N$-$M$-bimodule too
(see 3.1.11 in \cite{BLM}), and similarly for $X$.
We claim that $Y$ is selfdual as a right module if and
only $Y$ is selfdual as a left module.    
Indeed, if $Y$ is selfdual as a right $M$-module, then $X$ and $Y$ are
dual operator spaces  by  Lemma \ref{7prez}.
Hence using the left version of Corollary 
\ref{7masd},  we see that $Y$ is selfdual 
 as a left $N$-module.  Since $N$ is the 
multiplier algebra of the appropriate ideal,
it follows from 
Lemma \ref{wsid} that $Y$ is $w^*$-full as a left $N$-module.   

\begin{proposition} \label{newp}  Let 
$X$ and $Y$ be as in the last paragraph.  Then  $X$ and 
$Y$ have operator space preduals and are normal  
dual operator bimodules over $M$ and $N$, if and only if
$Y$ is selfdual and its canonical $w^*$-topology
as a selfdual right module  (mentioned after Lemma \ref{7prez}),
  agrees with its canonical $w^*$-topology as a
selfdual left module, and similarly for $X$.
\end{proposition}

\begin{proof}  ($\Leftarrow$) \ Follows from 
Lemma \ref{7prez}, and the method of
proof of (5) of that result.

($\Rightarrow$) \  Assuming
$X$ and
$Y$ have operator space preduals, we will refer to
the associated $w^*$-topologies as
{\em the original $w^*$-topologies} of $X$ and
$Y$.   By Theorem \ref{7osz}, $Y$ is
selfdual as a right $M$-module.   To say that a bounded 
net $(y_t)$ converges to $y \in Y$ in the $w^*$-topology 
mentioned after Lemma \ref{7prez},
is to say that $(x',y_t) \to (x',y)$ in the $w^*$-topology
of $M$, for all $x' \in X$.  By Corollary \ref{fromm}, this implies that
\begin{equation} \label{wneq}
y' (x',y_t) \to y' (x',y) \;   \textrm{in the original
$w^*$-topology of} \; Y , \;  \textrm{for all} \; x' \in X, y' \in Y .
\end{equation}
In fact it is equivalent to (\ref{wneq}), since if 
(\ref{wneq}) holds, and if $((x',y_{t_\mu}))_\mu$ 
is  a $w^*$-convergent
subnet of $((x',y_t))_t$ with limit $m \in M$,
then by Corollary \ref{fromm}, $y' (x',y_{t_\mu})$ 
converges to $y' m$.
This implies that $y' m = y' (x',y)$ for all
$y' \in Y$, so that $m = (x',y)$.  Hence 
$(x',y_t) \to (x',y)$ in the $w^*$-topology of $M$.
By Corollary \ref{fromm},
 if $y_t \to y$ in the original $w^*$-topology of
$Y$, then (\ref{wneq}) holds.  Conversely, if 
(\ref{wneq}) holds, then $y_t \to y$ in the original $w^*$-topology,
by a $w^*$-convergent subnet argument similar to the one 
we just used above.
(For if a subnet of $(y_t)$ converged with limit $y''$, say,
then using Corollary \ref{fromm} as above shows that
$y' (x',y'') = y' (x',y)$ for all such $x', y'$.
This implies that $y'' = y$.)

We have now shown that the canonical 
$w^*$-topology (mentioned after Lemma \ref{7prez})
of $Y$ as a right module agrees with its original $w^*$-topology.
By a symmetrical argument, this agrees with the canonical
$w^*$-topology as a left module.  Similarly, for $X$.
  \end{proof}   

The equivalent conditions in the last result are automatic in 
the $W^*$-algebra case, but not more generally.
If these conditions are satisfied, then we call $Y$ 
a {\em  tight $w^*$-equivalence $N$-$M$-bimodule},
and we say that $M$ and $N$ are {\em tightly
Morita $w^*$-equivalent}.  It then follows as in 
the second paragraph of the proof of Theorem \ref{7osz},
using also Lemma \ref{wsid},
 that $N \cong CB_M(Y)$ completely isometrically.
The isomorphism here takes $n \in N$ to the map $y \mapsto
ny$ on $Y$.  It is easy to argue, as in Lemma \ref{7prez} (5),
that this isomorphism is $w^*$-continuous.
Hence by the Krein-Smulian theorem it is a $w^*$-homeomorphism.
Thus, just as in the selfadjoint theory,
we can forget about $N$, and instead work with 
$CB_M(Y)$ (which equals $\Bdb_M(Y)$), when convenient. 
 
Conversely, we have:

\begin{theorem} \label{onweq}
Let $Y$ be a selfdual right rigged-equivalence module
over a unital dual operator algebra $M$.
Then $Y$ is a left $w^*$-full selfdual $CB_M(Y)$-rigged 
module.  Also, $Y$ implements a tight Morita
$w^*$-equivalence between  $CB_M(Y)$ and 
${\mathcal I}^w$.   
In particular, if $Y$ is also a right $w^*$-full 
$M$-module then $Y$ implements a tight Morita
$w^*$-equivalence between  $CB_M(Y)$ and $M$.
\end{theorem}

\begin{proof}
By Lemma \ref{7prez}, $Y$ and $X$ are
dual operator spaces.
By Corollary \ref{7prez2}, we have  
 $CB_M(Y) = \Bdb_M(Y)$, and
this is a dual operator algebra.
As we said at the end of the proof of Theorem \ref{7osz},
this space also equals $M(\Kdb_M(Y))$.
  The `left-hand variant' of Theorem \ref{7osz} says that 
$Y$ is a selfdual left $CB_M(Y)$-rigged
module, and it is $w^*$-full by 
Lemma \ref{7jism}.  The other assertions follow immediately 
from the definition of tight
Morita $w^*$-equivalence, and Lemma \ref{7prez} (5).
\end{proof}

{\bf Examples.}  Examples  of tight  Morita
$w^*$-equivalence, and therefore of selfdual right 
rigged-equivalence modules,  are not hard to find.
We list just three, omitting details: 

\smallskip

(1) \   $W^*$-algebras are  Morita equivalent 
in the sense of \cite{MEFC}, if and only if they
are tightly Morita $w^*$-equivalent.  This follows from
the definition in e.g.\ 8.5.12 of \cite{BLM}, and the fact
that for $C^*$-algebras, Rieffel's notion of strong Morita 
equivalence coincides with the one in \cite{COOM}
(see Chapter 6 of that reference).

\smallskip

(2) \ If $A$ and $B$ are any two unital operator 
algebras which are strongly Morita
equivalent in the sense of  \cite{COOM}, then $A^{**}$ and $B^{**}$ are tightly
Morita $w^*$-equivalent.   We omit the proof, which uses the method of 
8.5.32 in \cite{BLM}.
 
\smallskip

(3) \ Let $\eta$ be a fixed vector in a Hilbert space $H$.
The set of bounded operators on $H$ which have $\eta$ as an eigenvector,
 is a unital dual algebra which is tightly Morita $w^*$-equivalent
to the upper triangular $2 \times 2$ matrices.
The associated  $w^*$-equivalence bimodules may be taken to be
the set of operators from $\Cdb^2$ to $H$ taking  
the vector $e_1$ to a scalar multiple of $\eta$,
and the set of operators from $H$ to $\Cdb^2$ taking $\eta$
to a scalar multiple of $e_1$.  
 
\medskip

We next show that any selfdual right rigged-equivalence module
$X$ over a unital dual operator algebra $N$,
occurs as a `corner' of a unital dual operator algebra
${\mathcal L}$.
Note that the `right $N$-rigged sum' $X \oplus_c N$ is a right 
$N$-rigged module,
which is clearly selfdual.  The conjugate
left $N$-rigged module is $Y \oplus_r N$,
where $Y = \tilde{X}$  (see \cite[Section 4]{AGOH}).  Therefore,
by Lemma \ref{7prez},
 $X \oplus_c N$ is a dual operator space, and it is easy to 
check using Lemma \ref{7prez} (3) that the 
containments of $X$ and $N$ in this latter space are
$w^*$-homeomorphisms.  Let $p$ be the projection from
$X \oplus_c N$ onto $X \oplus 0$.
 Thus ${\mathcal L} = M_l(X \oplus_c N) = CB_N(X \oplus_c N)$ is 
a dual operator algebra,
by  \cite[Corollary 3.2]{MADO}.  The four corners of
${\mathcal L}$ are $X , Y, N$, and $M \cong CB_N(X)$;
indeed $X = p {\mathcal L} (1-p)$. 
By  Lemma \ref{7prez} (2) and Corollary \ref{7prez2} (2),
 the $w^*$-topologies on $X$ and $Y$
inherited from ${\mathcal L}$
coincide with the original ones.   Similarly for the 
other corners.

It is now clear that one has a theory 
that is simultaneously the appropriate 
`$w^*$-topology version' of 
 much of the theory in \cite{COOM}, and a generalization of  
much of the $C^*$-algebraic theory of 
weak Morita equivalence and $W^*$-modules
(see \cite[Section 8.5]{BLM}).   
Moreover it is clear that operator space multipliers
play an important role 
in this theory.   Generalizing many 
of the other results in the selfadjoint variant of the 
theory, is now essentially a routine exercise.
For example, one may show,
 analogously to a result due to Rieffel
in the $W^*$-algebra case,
that any selfdual rigged-equivalence module over a unital dual
operator algebra $M$ is of the form
$_RB(K,H)$, for a suitable Hilbert module $K$ over $M$, and
a Hilbert $R$-module $H$, where $R$ is the commutant
of $M$ in $B(K)$.   The argument follows the lines of
that of 8.5.37 and 8.5.32 in \cite{BLM},  
but using also the double commutant theorem for nonselfadjoint
operator algebras of Blecher and Solel (e.g.\ see
3.2.14 in \cite{BLM}).
  We will not prove it here,
since this result would take us 
away from the main themes of the present paper.

The main obstacles to the nonselfadjoint variant
of weak Morita equivalence presented here, that we 
see at this point, are twofold.
First, it is not clear, and probably is not true in general,
that a dual unital operator algebra $M$ is
always tightly $w^*$-Morita equivalent to $M \bar{\otimes} B(H)$,
if $H$ is an infinite dimensional Hilbert space.   
This is because the space $Y = M \bar{\otimes} H^c$,
the `first column' of $M \bar{\otimes} B(H)$,
is not a rigged module over $M$, in general,
unlike the $W^*$-algebra case.   Presumably this latter
deficiency may be fixed by considering weaker forms
of the rigged module definition.  However this will not
really help: this very natural $M$-module $Y$ is not 
even selfdual---there may exist completely bounded
$M$-module maps from $Y$ to $M$ which are 
not given by `left multiplication with a row in $M \bar{\otimes} H^r$',
where the latter space is the `first row' of $M \bar{\otimes} B(H)$.
An example of such is easy to construct in the case that  $M$ is
 the subalgebra of 
$M_2(B(\ell^2))$ with
$0$ in the 2-1 entry, scalars
on the main diagonal, and an element from $B(\ell^2)$ in the 1-2 entry.
This shows that any decent theory of 
selfdual modules over nonselfadjoint algebras has to 
either exclude such examples, or replace 
completely bounded
$M$-module maps by $w^*$-continuous ones (which
somewhat defeats the point of `selfduality'), or perhaps
by multipliers in the sense of the second Part of \cite{MCAA}.    
The second obstacle is it seems not to be true
in full generality, that the second dual    
of a strong Morita equivalence $A$-$B$-bimodule
in the sense of \cite{COOM}, is a tight $w^*$-equivalence
$A^{**}$-$B^{**}$-bimodule in the sense above.

Some of these problems are easily resolvable,
at the expense of introducing other problems,
 if one instead uses a different approach to
$w^*$-Morita theory.  In fact there are several  
such alternative approaches.  First, one could vary the 
theory above by allowing the `special approximate
identities' found in the theory
to converge in the point-$w^*$ topology as opposed
to the point-norm  topology:
for example $\sum_{k = 1}^{n_\alpha} y_k^\alpha (x_k^\alpha , y)
\to y$  in the $w^*$ topology for all $y \in Y$.  
Second, another completely different approach is to 
base the entire theory
on a  (not yet developed)
nonselfadjoint dual operator algebra variant
of the Haagerup module tensor product
(cf.\  \cite{SOMA}).  However, both of
these approaches seems to present other, different,
problems.  For example, it seems certain that one cannot obtain,
by such approaches,
analogues of many of our results here. 
 
\medskip
 
{\bf Acknowledgments:}   We thank Christian Le Merdy and Vrej
 Zarikian for several conversations and inputs.

 \end{document}